\documentclass[12pt,reqno]{amsart}
\usepackage{amsmath,amssymb,amsfonts,amscd,latexsym,amsthm,mathrsfs,verbatim}
\theoremstyle{remark}

\renewcommand{\d}{{\mathrm d}}

\renewcommand{\Re}{\operatorname{Re}}
\begin{document}

\title{Lost in translation}

\author{Wadim Zudilin}
\address{School of Mathematical and Physical Sciences, The University of Newcastle, Callaghan, NSW 2308, Australia}
\email{wadim.zudilin@newcastle.edu.au}

\thanks{The author is supported by the Australian Research Council.} %grant DP110104419.

\date{26 September 2012}


\begin{abstract}
We explain the use and set grounds about applicability
of algebraic transformations of arithmetic hypergeometric series for proving
Ramanujan's formulae for $1/\pi$ and their generalisations.
\end{abstract}

\dedicatory{In memory of Herb Wilf}

\maketitle
%==================================================

The principal goal of this note is to set some grounds about applicability
of algebraic transformations of (arithmetic) hypergeometric series for proving
Ramanujan's formulae for $1/\pi$ and their numerous generalisations.
The technique was successfully used in quite different situations
\cite{CZ10,Rog09,WZ12,Zu07a,Zu07b} and was dubbed as `translation
method' by J.~Guillera, although the name does not give any clue
about the method itself.

\bigskip
Consider the following problem: \emph{Show that}
\begin{equation}
\sum_{n=0}^\infty\frac{(4n)!}{n!^4}
(3+40n)\cdot\frac{1}{28^{4n}}
=\frac{49}{3\sqrt3\pi}\,.
\label{wz01}
\end{equation}

\medskip
\noindent
\textsl{Step 0}. It comes as a useful rule: prior to any attempts to prove an
identity verify it numerically. The convergence of the series on the left-hand side
of~\eqref{wz01} is reasonably fast (more than 3 decimal places per term), so
you shortly convince yourself that the both sides are

{\small
$$
3.001679541740867825117222046370611403163548615329487998574326\dotsc.
$$
}\relax

\medskip
\noindent
\textsl{Step 1}. Series of the type given in \eqref{wz01} should be
quite special. With a little search you identify
\begin{equation}
\sum_{n=0}^\infty\frac{(4n)!}{n!^4}\biggl(\frac x{256}\biggr)^n
={}_3F_2\biggl(\begin{matrix} \frac14, \, \frac12, \, \frac34 \\ 1, \, 1 \end{matrix}\biggm| x \biggr)
=\sum_{n=0}^\infty\frac{\bigl(\frac14\bigr)_n\bigl(\frac12\bigr)_n\bigl(\frac34\bigr)_n}{(1)_n(1)_n}\,\frac{x^n}{n!}\,,
\label{wz02}
\end{equation}
a hypergeometric series, where the notation $(a)_n$ (\emph{Pochhammer's symbol} or \emph{shifted factorial})
stands for $\Gamma(a+n)/\Gamma(a)=a(a+1)\dotsb(a+n-1)$. A generalised hypergeometric series
\begin{equation*}
{}_mF_{m-1}\biggl(\begin{matrix}
a_1, \, a_2, \, \dots, \, a_m \\
b_2, \, \dots, \, b_m \end{matrix}\biggm|x\biggr)
:=\sum_{n=0}^\infty
\frac{(a_1)_n(a_2)_n\dotsb(a_m)_n}
{(b_2)_n\dotsb(b_m)_n}\,
\frac{x^n}{n!}
\end{equation*}
is an object of intensive study since Euler \cite{Bai35,Sl66}; one of its important
properties is the linear differential equation
\begin{equation}
\biggl(\biggl(x\frac{\d}{\d x}\biggr)\prod_{j=2}^m\biggl(x\frac{\d}{\d x}+b_j-1\biggr)
-x\prod_{j=1}^m\biggl(x\frac{\d}{\d x}+a_j\biggr)\biggr)F=0
\label{wz03}
\end{equation}
satisfied by the series. The required identity~\eqref{wz01} can be therefore transformed
to the more conceptual form
\begin{equation}
\sum_{n=0}^\infty\frac{\bigl(\frac14\bigr)_n\bigl(\frac12\bigr)_n\bigl(\frac34\bigr)_n}{n!^3}
\,\frac{3+40n}{7^{4n}}
=\biggl(3+40x\frac{\d}{\d x}\biggr)
{}_3F_2\biggl(\begin{matrix} \frac14, \, \frac12, \, \frac34 \\ 1, \, 1 \end{matrix}\biggm| x \biggr)\bigg|_{x=1/7^4}
=\frac{49}{3\sqrt3\pi}\,.
\label{wz04}
\end{equation}

\medskip
\noindent
\textsl{Step 2}. Convince yourself that identities of the wanted type
are known in the literature. In fact, they are known
for almost a century after Ramanujan's publication~\cite{Ra14};
identity \eqref{wz01} is equation~(42) there.
Ramanujan did not indicate how he arrived at his series
but left some hints that these series belong to what is now known as
`the theories of elliptic functions to alternative bases'.
The first proofs of Ramanujan's identities and
their generalisations were given by the Borweins~\cite{BB87} and
Chudnovskys~\cite{CC88}. Those proofs are however too lengthy to be included here.
Note that Ramanujan's list in~\cite{Ra14} does not include
the slowly convergent example
\begin{equation}
\sum_{n=0}^\infty\frac{\bigl(\frac12\bigr)_n^3}{n!^3}
(1+4n)\,(-1)^n
=\biggl(1+4x\frac{\d}{\d x}\biggr)
{}_3F_2\biggl(\begin{matrix} \frac12, \, \frac12, \, \frac12 \\ 1, \, 1 \end{matrix}\biggm| x \biggr)\bigg|_{x=-1}
=\frac2\pi\,,
\label{wz05}
\end{equation}
which was shown to be true by G.~Bauer~\cite{Ba59} already in~1859.
Bauer's proof makes no reference to sophisticated theories and is much
shorter, although does not seem to be generalisable to the other entries
from~\cite{Ra14}. In fact, D.~Zeilberger assisted by his automatic collaborator
S.\,B.~Ekhad \cite{EZ94} came up in 1994 with a short proof of~\eqref{wz05}
verifiable by a computer. The key is a use of a simple telescoping argument
(this part is completely automated by the great Wilf--Zeilberger (WZ) machinery~\cite{PWZ97})
and an advanced theorem due to Carlson~\cite[Chap.~V]{Bai35}; the proof is reproduced in~\cite{Zu08}.
Quite recently, J.~Guillera advocated \cite{Gu02,Gu06,Gu10,Gu11} the method
from~\cite{EZ94} and significantly extended the outcomes;
he showed, for example, that many other Ramanujan's identities for $1/\pi$
can be proven completely automatically. Note however that
\eqref{wz01} is one of `WZ resistent' identities. To overcome this technical
difficulty, below we reduce the identity to the simpler one~\eqref{wz05}.

\medskip
\noindent
\textsl{Step 3}.
Use your favourite computer algebra system (CAS) to verify
the hypergeometric identity
\begin{equation}
{}_3F_2\biggl(\begin{matrix} \frac12, \, \frac12, \, \frac12 \\ 1, \, 1 \end{matrix}\biggm| x \biggr)
=r\cdot{}_3F_2\biggl(\begin{matrix} \frac14, \, \frac12, \, \frac34 \\ 1, \, 1 \end{matrix}\biggm| y \biggr)
\label{wz06}
\end{equation}
where $y=y(x)=-\frac1{1024}x^3+O(x^4)$ and $r=r(x)=1+\frac18x+\frac{27}{512}x^2+O(x^3)$
are algebraic functions determined by the equations
{\small
\begin{align*}
&
(x^2-194x+1)^4y^4
\\ &\quad
+16(4833x^6+2029050x^5+47902255x^4-92794388x^3
\\ &\quad\qquad
+47902255x^2+2029050x+4833)xy^3
\\ &\quad
-96(3328x^6-623745x^5+3837060x^4-6470150x^3
\\ &\quad\qquad
+3837060x^2-623745x+3328)xy^2
\\ &\quad
+256(1024x^6-1152x^5+225x^4-2x^3+225x^2-1152x+1024)xy+256x^4
=0
\end{align*}
}\relax
and
\begin{align*}
&
(x^2-194x+1)^2r^8+4(61x^2+25798x+61)(x-1)r^6
\\ &\quad
+486(41x^2-658x+41)r^4+551124(x-1)r^2+531441
=0.
\end{align*}
To do this you (and your CAS) are expected to use the linear differential
equations~\eqref{wz03} for the involved hypergeometric functions and
generate any-order derivatives of $y$ and $r$ with respect to~$x$ by appealing
to the implicit functional equations. To summarise, you have to check that
the both sides of \eqref{wz06} satisfy the same (3rd order) linear
differential equation in~$x$ with algebraic function coefficients and
then compare the first few coefficients in the expansions in powers of~$x$.
Note that $x=-1$ corresponds to $y=1/7^4$ (cf.~\eqref{wz05} vs.~\eqref{wz04}),
and this is the reason behind considering the sophisticated functional identity~\eqref{wz06}.

The task on this step does not look humanly pleasant, and there is a
(casual) trick to verify~\eqref{wz06} by parameterising $x$, $y$ and~$r$:
\begin{gather*}
x=-\frac{4p(1-p)(1+p)^3(2-p)^3}{(1-2p)^6},
\quad
y=\frac{16p^3(1-p)^3(1+p)(2-p)(1-2p)^2}{(1-2p+4p^3-2p^4)^4},
\\
r=\frac{(1-2p)^3}{1-2p+4p^3-2p^4}.
\end{gather*}
Choosing $p=(1-\sqrt{45-18\sqrt6})/2$ we obtain $x=-1$ and $y=1/7^4$.
(The modular reasons behind this parametrisation can be found in
\cite[Lemma 5.5 on p.~111]{Be98} where the $p$ there is the negative
of our~$p$.)

\medskip
\noindent
\textsl{Step 4}. By differentiating identity~\eqref{wz06} with respect to $x$
and combining the result with~\eqref{wz06} itself we see that
\begin{equation}
\biggl(a+bx\frac{\d}{\d x}\biggr)
{}_3F_2\biggl(\begin{matrix} \frac12, \, \frac12, \, \frac12 \\ 1, \, 1 \end{matrix}\biggm| x \biggr)
=\biggl(a+bx\frac{\d r}{\d x}+b\frac{rx}y\,\frac{\d y}{\d x}\cdot y\frac{\d}{\d y}\biggr)
\cdot{}_3F_2\biggl(\begin{matrix} \frac14, \, \frac12, \, \frac34 \\ 1, \, 1 \end{matrix}\biggm| y \biggr);
\label{wz07}
\end{equation}
again, the derivatives $\d y/\d x$ and $\d r/\d x$ are read
from the implicit functional equations. An alternative (but simpler) way is using
the parametrisations $x(p)$, $y(p)$ and $r(p)$. Taking $a=1$, $b=4$ and $x=-1$
in~\eqref{wz07} you recognise the familiar Bauer's (WZ easy) identity~\eqref{wz05}
on the left-hand side; the right-hand side counterpart is nothing but~\eqref{wz04}.

\bigskip
\noindent
\textsl{Comments}.
The story exposed above is general enough to be used in other situations for
proving \emph{some} other formulae for $1/\pi$. The setup can be as follows.
Assume we already have an identity
$$
\biggl(a+bx\frac{\d}{\d x}\biggr)F(x)\bigg|_{x=x_0}=\mu,
$$
where $a$, $b$, $x_0$ and $\mu$ are certain (simple or at least
arithmetically significant) numbers, and
$F(x)$~is an (arithmetic) series. Furthermore, assume we have
a transformation $F(x)=rG(y)$ with $r=r(x)$ and $y=y(x)$ differentiable
at $x=x_0$. Then
$$
\biggl(\hat a+\hat by\frac{\d}{\d y}\biggr)G(y)\bigg|_{y=y_0}=\mu,
$$
where
$$
\hat a=a+bx\frac{\d r}{\d x}\bigg|_{x=x_0}, \quad
\hat b=b\frac{rx}y\,\frac{\d y}{\d x}\bigg|_{x=x_0}, \quad\text{and}\quad
y=y_0.
$$
There is, of course, no magic in this result: it is just the standard `chain rule'.

The applicability of this simple argument heavily rests on existence of
transformations like \eqref{wz06}. This in turn is based on the modular
origin~\cite{BB87,CC12,CC88,Zu08} of Ramanujan's identities for $1/\pi$: any
such identity can be written in the form
\begin{equation}
\biggl(a+bx\frac{\d}{\d x}\biggr)F(x)\bigg|_{x=x_0}=\frac c\pi,
\qquad a,b,c,x_0\in\mathbb Q,
\label{wz08}
\end{equation}
where $F(x)$ is an \emph{arithmetic hypergeometric series}~\cite{Zu11} satisfying
a 3rd order linear differential equation. In other words, for a certain
modular function $x=x(\tau)$ (not uniquely defined!) the function $F(x(\tau))$
is a modular form of weight~2.
The theory of modular forms provides us with the knowledge that any two modular
forms are algebraically dependent; thus, whenever we have another arithmetic
hypergeometric series $G(y)$ and a related modular parametrisation $y=y(\tau)$,
the modular functions $y(\tau)$ and $G(y(\tau))/F(x(\tau))$ are algebraic
over $\mathbb Q[x(\tau)]$. Another warrants of the theory is an algebraic
dependence over $\mathbb Q$ of $x(\tau)$ and $x((A\tau+B)/(C\tau+D))$ for any
$\bigl(\begin{smallmatrix} A & B \\ C & D\end{smallmatrix}\bigr)\in SL_2(\mathbb Q)$.
On the other hand, there is no other source known for such algebraic dependency;
the functions $x(\tau)$ and $x(A\tau)$, $A>0$, are algebraically dependent if
and only if $A$~is rational.

The above arithmetic constraints impose the natural restriction on
$\tau_0$ from the upper half-plane $\Re\tau>0$ to satisfy $x(\tau_0)=x_0$ in~\eqref{wz08}.
Namely, $\tau_0$ is an (imaginary) quadratic irrationality, $\tau_0\in\mathbb Q[\sqrt{-d}]$
for some positive integer~$d$. But then $(A\tau_0+B)/(C\tau_0+D)$ belongs to the same
quadratic extension of~$\mathbb Q$ for any
$\bigl(\begin{smallmatrix} A & B \\ C & D\end{smallmatrix}\bigr)\in SL_2(\mathbb Q)$,
so whatever transformation $F(x)=rG(y)$ (of modular origin) we use, the modular
arguments of $x(\tau)$ and $y(\tau)$ have to be tied by an $SL_2(\mathbb Q)$ linear-fractional
transform. In the examples~\eqref{wz04} and~\eqref{wz05} we have both arguments
belonging to $\mathbb Q[\sqrt{-2}]$, therefore an algebraic transformation must exist,
and this is confirmed by~\eqref{wz06} mapping the corresponding $x(\tau_0)=-1$
into $y(3\tau_0)=1/7^4$ where $\tau_0=(1+\sqrt{-2})/2$. There is however no way known
to `translate' identities~\eqref{wz04} and~\eqref{wz05} to either
\begin{align*}
\sum_{n=0}^\infty\frac{\bigl(\frac12\bigr)_n^3}{n!^3}
(1+6n)\,\frac1{4^n}
&=\frac4\pi
\\ \intertext{or}
\sum_{n=0}^\infty\frac{\bigl(\frac16\bigr)_n\bigl(\frac12\bigr)_n\bigl(\frac56\bigr)_n}{n!^3}
(13591409+545140134n)\cdot\frac{(-1)^n}{53360^{3n+2}}
&=\frac3{2\sqrt{10005}\pi}\,,
\end{align*}
as the corresponding modular arguments lie in the fields $\mathbb Q[\sqrt{-3}]$ and
$\mathbb Q[\sqrt{-163}]$, respectively. We refer the interested reader to~\cite{CC12}
for exhausting lists of `rational' (in the sense of~$x_0$) identities which express $1/\pi$
by means of general hypergeometric-type series; the details of the modular machinery
are greatly explained there.

In a sense, to make the `translation method' work we first should carefully examine
the underlying modular parametrisations. On the other hand, there are situations
when we know (or can produce~\cite{ASZ11}) the algebraic transformations without
having modularity at all. These are particularly useful in the context of
similar formulae for $1/\pi^2$ recently discovered by Guillera~\cite{Gu02,Gu06,Gu11}.

There is a $p$-adic counterpart of the Ramanujan-type identities for $1/\pi$
and $1/\pi^2$ which we review in~\cite{Zu09}. It seems likely that
the algebraic transformation machinery is generalisable to those situations
as well but, for the moment, no single example of this is known.

\medskip
\textbf{Acknowledgements.} I would like to thank Shaun Cooper for his useful
suggestions which helped me to improve on an earlier draft of this note.

%==================================================

\end{document}